\documentclass{amsart}

\usepackage{amsfonts,amssymb,amscd,amstext}
\usepackage[a4paper,hmargin=3.6cm,vmargin=5cm]{geometry}
\usepackage{mathpazo}

\usepackage{graphicx}
\usepackage{color}

\renewcommand{\le}{\leqslant}
\renewcommand{\ge}{\geqslant}
\newcommand{\ptl}{\partial}

\newcommand{\rr}{{\mathbb{R}}}

\newcommand{\la}{\lambda}

\newcommand{\hh}{{\mathbb{H}}}

\newcommand{\sph}{{\mathbb{S}}}
\newcommand{\zz}{\mathbb{Z}}
\newcommand{\pp}{\mathcal{P}}
\newcommand{\h}{\mathcal{H}}

\newcommand{\sub}{\subset}
\newcommand{\subeq}{\subseteq}

\newcommand{\escpr}[1]{\big<#1\big>}
\newcommand{\Sg}{\Sigma}
\newcommand{\Om}{\Omega}
\newcommand{\al}{\alpha}
\newcommand{\eps}{\varepsilon}

\newcommand{\ga}{\gamma}
\newcommand{\Ga}{\Gamma}

\newcommand{\wt}{\widetilde}

\DeclareMathOperator{\divv}{div}

\DeclareMathOperator{\intt}{int}

\DeclareMathOperator{\sgn}{sgn}

\newtheorem{theorem}{Theorem}[section]
\newtheorem{proposition}[theorem]{Proposition}
\newtheorem{lemma}[theorem]{Lemma}

\theoremstyle{definition}

\newtheorem{example}[theorem]{Example}

\theoremstyle{remark}

\numberwithin{equation}{section}

\begin{document}

\title[Examples of area-minimizing surfaces]{Examples of
area-minimizing surfaces in the subriemannian Heisenberg group $\hh^1$
with low regularity}

\author[M.~Ritor\'e]{Manuel Ritor\'e} \address{Departamento de
Geometr\'{\i}a y Topolog\'{\i}a \\
Universidad de Granada \\ E--18071 Granada \\ Espa\~na}
\email{ritore@ugr.es}

\date{\today}

\thanks{Research supported by MEC-Feder grant MTM2007-61919}

\begin{abstract}
We give new examples of entire area-minimizing $t$-graphs in the
subriemannian Heisenberg group $\hh^1$.  Most of the examples are
locally lipschitz in Euclidean sense.  Some regular examples have
prescribed singular set consisting of either a horizontal line or a
finite number of horizontal halflines extending from a given point.
Amongst them, a large family of area-minimizing cones is obtained.
\end{abstract}

\maketitle

\thispagestyle{empty}

\section{Introduction}

Variational problems related to the subriemannian area in the
Heisenberg group $\hh^1$ have received great attention recently.  A
major question in this theory is the regularity of minimizers.  A
related one is the construction of examples with low regularity
properties.  The study of minimal surfaces in subriemannian geometry
was initiated in the paper by Garofalo and Nhieu \cite{gn}.  Later
Pauls \cite{pauls-examples} constructed minimal surfaces in $\hh^1$ as
limits of minimal surfaces in Nil manifolds, the riemannian Heisenberg
groups.  Cheng, Hwang and Yang \cite{chy1} have studied the weak
solutions of the minimal surface equation for $t$-graphs and have
proven existence and uniqueness results.  Regularity of minimal surfaces,
assuming that they are least $C^1$, has been treated in the papers by
Pauls \cite{pauls-regularity} and Cheng, Hwang and Yang \cite{chy2}.
We would like also to mention the recently distributed notes by
Bigolin and Serra Cassano \cite{bsc}, where they obtain regularity
properties of an $\hh$-regular surface from regularity properties
of its horizontal unit normal.  Interesting examples of minimal
surfaces which are not area-minimizing are obtained in
\cite{danielli-2006b}.  See also \cite{danielli-2006}. Smoothness of
lipschitz minimal intrinsic graphs in Heisenberg groups $\hh^n$, for
$n>1$, has been recently obtained by Capogna, Citti and Manfredini
\cite{ccm}.

Characterization in $\hh^1$ of solutions of the Bernstein problem for
$C^2$ surfaces has been obtained by Cheng, Hwang, Malchiodi and Yang
\cite{chmy}, and Ritor\'e and Rosales \cite{rr2} for $t$-graphs, and
by Barone Adessi, Serra Cassano and Vittone \cite{bscv} and Garofalo
and Pauls \cite{gp} for vertical graphs.

Additional contributions concerning variational problems related to
the subriemannian area in the Heisenberg groups include \cite{pansu2},
\cite{ascv}, \cite{chmy}, \cite{chy1}, \cite{chy2}, \cite{MR2313532},
\cite{MR2313532}, \cite{MR2287539}, \cite{MR2032504},
\cite{MR1984849}, \cite{fssc}, \cite{MR1711057}, \cite{msc},
\cite{rr2}.  The recent monograph by Capogna, Danielli, Pauls and
Tyson \cite{cdpt} gives a recent overview of the subject with an
exhaustive list of references.  We would like to stress that, in
$\hh^1$, the condition $H\equiv 0$ is not enough to guarantee that a
given surface of class $C^2$ is even a stationary point for the area
functional, see Ritor\'e and Rosales \cite{rr2}, and Cheng, Hwang and
Yang \cite{chy1} for minimizing $t$-graphs.

The aim of this paper is to provide new examples in $\hh^1$ of
Euclidean locally lipschitz area-minimizing entire graphs over the
$xy$-plane.

In section~\ref{sec:line} we construct the basic examples.  We start
from a given horizontal line $L$, and a monotone angle function
$\alpha:L\to (0,\pi)$ over this line.  For each $p\in L$, we consider
the two horizontal halflines extending from $p$ making an angle
$\pm\alpha(p)$ with $L$.  We prove that in this way we always obtain
an entire graph over the $xy$-plane which is Euclidean locally
lipschitz and area-minimizing.  The angle function $\alpha$ is only
assumed to be continuous and monotone.  Of course, further regularity
on $\alpha$ yields more regularity on the graph.  In case $\alpha$ is
at least $C^2$ we get that the associated surface is $C^{1,1}$.

The surfaces in section~\ref{sec:line} are the building blocks for our
next construction in section~\ref{sec:halflines}.  We fix a point
$p\in\hh^1$, and a family of counter-clockwise oriented horizontal
halflines $R_{1}$, $\ldots$, $R_{n}$ extending from $p$.  We choose
the bisector $L_{i}$ of the wedge determined by $R_{i-1}$ and $R_{i}$,
and we consider angle functions $\alpha_{i}:L_{i}\to (0,\pi)$ which
are continuous, nonincreasing as a function of the distance to $p$,
and such that $\alpha(p)$ is equal to the angle between $L_{i}$ and
$R_{i}$.  For every $q\in L_{i}$, we consider the halflines extending
from $q$ with angles $\pm\alpha_{i}(q)$.  In this way we also a family
of area-minimizing $t$-graphs which are Euclidean locally lipschitz.
In case the obtained surface is regular enough we have that the
singular set is precisely $\bigcup_{i=1}^n L_{i}$.  If the angle
functions $\alpha_{i}$ are constant, then we obtain area-minimizing
cones (the original motivation of this paper), which are Euclidean
locally $C^{1,1}$ minimizers, and $C^{\infty}$ outside the singular
set $\bigcup_{i=1}^n L_{i}$.  For a single halfline $L$ extending from
the origin and an angle function $\alpha:L\to (0,\pi)$, continuous and
nonincreasing as a function of the distance to $0$, we patch the graph
obtained over a wedge of the $xy$-plane with the plane $t=0$ along the
halflines extending from $0$ making an angle $\alpha(0)$ with $L$.
When $\alpha$ is constant we get again an area-minimizing cone which
is Euclidean locally lipschitz.  These cones are a generalization of
the one obtained by Cheng, Hwang and Yang \cite[Ex.~7.2]{chy1}.

An interesting consequence of this construction is that we get a large
number of Euclidean locally $C^{1,1}$ area-minimizing cones with
prescribed singular set consisting on either a horizontal line or a
finite number of horizontal halflines extending from a given point.
It is an open question to decide if these examples are the only
area-minimizing cones, together with vertical halfspaces and the
example by Cheng, Hwang and Yang \cite[Ex.~7.2]{chy1} with a singular
halfline and its generalizations in the last section.  The importance
of tangent cones has been recently stressed in \cite{ambrosio-2008}.

\section{Preliminaries}

The \emph{Heisenberg group} $\hh^1$ is the Lie group $(\rr^3,*)$,
where the product $*$ is defined, for any pair of points $[z,t]$,
$[z',t']\in\rr^3\equiv\mathbb{C}\times\rr$, as
\[
[z,t]*[z',t']:=[z+z',t+t'+\text{Im}(z\overline{z}')], \qquad (z=x+iy).
\]
For $p\in\hh^1$, the \emph{left translation} by $p$ is the
diffeomorphism $L_p(q)=p*q$.  A basis of left invariant vector fields
(i.e., invariant by any left translation) is given by
\begin{equation*}
X:=\frac{\ptl}{\ptl x}+y\,\frac{\ptl}{\ptl t}, \qquad
Y:=\frac{\ptl}{\ptl y}-x\,\frac{\ptl}{\ptl t}, \qquad
T:=\frac{\ptl}{\ptl t}.
\end{equation*}
The \emph{horizontal distribution} $\mathcal{H}$ in $\hh^1$ is the
smooth planar one generated by $X$ and $Y$.  The \emph{horizontal
projection} of a vector $U$ onto $\mathcal{H}$ will be denoted by
$U_{H}$.  A vector field $U$ is called \emph{horizontal} if $U=U_H$.
A \emph{horizontal curve} is a $C^1$ curve whose tangent vector lies
in the horizontal distribution.

We denote by $[U,V]$ the Lie bracket of two $C^1$ vector fields $U$, $V$
on $\hh^1$. Note that $[X,T]=[Y,T]=0$, while
$[X,Y]=-2T$.  The last equality implies that $\mathcal{H}$ is a
bracket generating distribution. Moreover, by Frobenius Theorem we
have that $\mathcal{H}$ is nonintegrable.  The vector fields $X$ and $Y$
generate the kernel of the (contact) $1$-form
$\omega:=-y\,dx+x\,dy+dt$.

We shall consider on $\hh^1$ the (left invariant) Riemannian metric
$g=\escpr{\cdot\,,\cdot}$ so that $\{X,Y,T\}$ is an orthonormal basis
at every point, and the associated Levi-Civit\'a connection $D$.  The
modulus of a vector field $U$ will be denoted by $|U|$.

Let $\ga:I\to\hh^1$ be a piecewise $C^1$ curve defined on a compact
interval $I\sub\rr$.  The \emph{length} of $\ga$ is the usual
Riemannian length $L(\ga):=\int_{I}|\dot{\ga}|$, where $\dot{\ga}$ is
the tangent vector of $\ga$.  For two given points in $\hh^1$ we can
find, by Chow's connectivity Theorem \cite[p.  95] {Gr}, a horizontal
curve joining these points.  The \emph{Carnot-Carath\'edory distance}
$d_{cc}$ between two points in $\hh^1$ is defined as the infimum of
the length of horizontal curves joining the given points. A
\emph{geodesic} $\ga:\hh^1\to\rr$ is a horizontal curve which is a
critical point of length under variations by horizontal curves. They
satisfy the equation
\begin{equation}
\label{eq:geodesics}
D_{\dot{\ga}}\dot{\ga}+2\la\,J(\dot{\ga})=0,
\end{equation}
where $\la\in\rr$ is the \emph{curvature} of the geodesic, and $J$ is
the $\pi/2$-degrees oriented rotation in the horizontal distribution.
Geodesics in $\hh^1$ with $\la=0$ are horizontal straight lines. The
reader is referred to the section on geodesics in \cite{rr2} for
further details.

The \emph{volume} $|\Om|$ of a Borel set $\Om\subeq\hh^1$ is the
Riemannian volume of the left invariant metric $g$, which coincides
with the Lebesgue measure in $\rr^3$.  We shall denote this volume
element by $dv_{g}$.  The \emph{perimeter} of $E\subset\hh^1$ in an
open subset $\Om\subset\hh^1$ is defined as
\begin{equation}
\label{eq:pereom}
|\ptl E|(\Om):=\sup\bigg\{\int_{\Om}\divv U\,dv_{g} : U\ 
\text{horizontal and }C^1, |U|\le 1, \text{supp}(U)\subset\Om\bigg\},
\end{equation}
where $\text{supp}(U)$ is the support of $U$.  A set $E\subset\hh^1$
is of \emph{locally finite perimeter} if $\pp(E,\Om)<+\infty$ for any
bounded open set $\Om\subset\hh^1$.  A set of locally finite perimeter
has a measurable \emph{horizontal unit normal} $\nu_{E}$, that
satisfies the following divergence theorem \cite[Corollary~7.6]{fssc}:
if $U$ is a horizontal vector field with compact support, then
\[
\int_{E} \divv U\,dv_{g}=\int_{\hh^1} \escpr{U,\nu_{E}}\,d|\ptl E|.
\]
If $E\subset\hh^1$ has Euclidean lipschitz boundary, then
\cite[Corollary~7.7]{fssc}
\begin{equation}
\label{eq:perlip}
|\ptl E|(\Om)=\int_{\ptl E\cap\Om}|N_{H}|\,d\mathcal{H}^2,
\end{equation}
where $N$ is the outer unit normal to $\ptl E$, defined
$\mathcal{H}^2$-almost everywhere.  Here $\mathcal{H}^2$ is the
$2$-dimensional riemannian Hausdorff measure. 

Let $\Om\subset\hh^1$ be an open set. We say that $E\subset\hh^1$ of
locally finite perimeter is \emph{area-minimizing} in $\Om$ if, for
any set $F$ such that $E=F$ outside $\Om$ we have
\[
|\ptl E|(\Om)\le |\ptl F|(\Om).
\]
The following extension of the divergence theorem will be needed to
prove the area-minimizing property of sets of locally finite perimeter
\begin{theorem}
Let $E\subset\hh^1$ be a set of locally finite perimeter,
$B\subset\hh^1$ a set with piecewise smooth boundary, and $U$ a $C^1$
horizontal vector field in $\intt(B)$ that extends continuously to the
boundary of $B$. Then
\begin{equation}
\label{eq:div}
\int_{E\cap B}\divv U\,dv_{g}=\int_{B}\escpr{U,\nu_{E}}\,d|\ptl E|
+\int_{E} \escpr{U,\nu_{B}}\,d|\ptl B|.
\end{equation}
\end{theorem}

\begin{proof}
The proof is modelled on \cite[\S~5.7]{eg}. Let $s$ denote the
riemannian distance function to $\hh^1-B$. For $\eps>0$, define
\[
h_{\eps}(p):=
\begin{cases}
1, &\eps\le s(p), \\
s(p)/\eps, & 0\le s(p)\le \eps,
\end{cases}
\]
Then $h_{\eps}$ is a lipschitz function (in riemannian sense).  For any
smooth $h$ with compact support in $B$ we have
$\divv(hU)=h\divv(U)+\escpr{\nabla h,U}$.  By applying the divergence
theorem for sets of locally finite perimeter \cite{fssc} we get
\[
\int_{\hh^1} h\,\escpr{U,\nu_{E}}\,d|\ptl E|=\int_{E}h\,\divv(U)
+\int_{E} \escpr{\nabla h,U}.
\]
By approximation, this formula is also valid for $h_{\eps}$. Taking
limits when $\eps\to 0$ we have $\h_{\eps}\to\chi_{B}$. By the coarea 
formula for lipschitz functions
\[
\frac{1}{\eps}\,\int_{\{0\le s\le\eps\}} \chi_{E}\,\escpr{\nabla s,U}=
\frac{1}{\eps}\int_{0}^\eps\bigg\{\int_{\{s=r\}}\chi_{E}\,\escpr{\nabla
s,U}\,d\mathcal{H}^2\bigg\}\,dr,
\]
and, taking again limits when $\eps\to 0$ and calling $N_{B}$ to the
riemannian outer unit normal to $\ptl B$ (defined except on a small
set), we have
\[
\lim_{\eps\to 0}\int_{E}\escpr{\nabla h_{\eps},U}=\int_{\ptl B}\chi_{E}
\escpr{N_{B},U}\,d\mathcal{H}^2=\int_{E}\escpr{\nu_{B},U}\,d|\ptl B|.
\]
Hence \eqref{eq:div} is proved.
\end{proof}

For a $C^1$ surface $\Sg\sub\hh^1$ the \emph{singular set} $\Sg_0$
consists of those points $p\in\Sg$ for which the tangent plane
$T_p\Sg$ coincides with the horizontal distribution.  As $\Sg_0$ is
closed and has empty interior in $\Sg$, the \emph{regular set}
$\Sg-\Sg_0$ of $\Sg$ is open and dense in $\Sg$.  It was proved in
\cite[Lemme 1]{d2}, see also \cite[Theorem~1.2]{balogh}, that, for a
$C^2$ surface, the Hausdorff dimension with respect to the Riemannian
distance on $\hh^1$ of $\Sg_{0}$ is less than two.

If $\Sg$ is a $C^1$ oriented surface with unit normal vector $N$, then
we can describe the singular set $\Sg_0\sub\Sg$, in terms of $N_H$, as
$\Sg_{0}=\{p\in\Sg:N_H(p)=0\}$.  In the regular part $\Sg-\Sg_0$, we
can define the \emph{horizontal unit normal vector} $\nu_H$, as in
\cite{dgn}, \cite{rr1} and \cite{gp} by
\begin{equation}
\label{eq:nuh}
\nu_H:=\frac{N_H}{|N_H|}.
\end{equation}
Consider the \emph{characteristic vector field} $Z$ on $\Sg-\Sg_0$
given by
\begin{equation}
\label{eq:zeta}
Z:=J(\nu_H).
\end{equation}
As $Z$ is horizontal and orthogonal to $\nu_H$, we conclude that $Z$
is tangent to $\Sg$.  Hence $Z_{p}$ generates the intersection of
$T_{p}\Sg$ with the horizontal distribution.  The integral curves of
$Z$ in $\Sg-\Sg_0$ will be called \emph{characteristic curves} of
$\Sg$.  They are both tangent to $\Sg$ and horizontal. Note that these
curves depend on the unit normal $N$ to $\Sg$.   If we define
\begin{equation}
\label{eq:ese}
S:=\escpr{N,T}\,\nu_H-|N_H|\,T,
\end{equation}
then $\{Z_{p},S_{p}\}$ is an orthonormal basis of $T_p\Sg$ whenever
$p\in\Sg-\Sg_0$.

In the Heisenberg group $\hh^1$ there is a one-parameter group of
\emph{dilations} $\{\varphi_s\}_{s\in\rr}$ generated by the vector
field
\begin{equation}
\label{eq:w}
W:=xX+yY+2tT.
\end{equation}
We may compute $\varphi_{s}$ in coordinates
to obtain
\begin{equation}
\label{eq:dilations}
\varphi_{s}(x_{0},y_{0},t_{0})=(e^sx_{0},e^sy_{0},e^{2s}t_{0}).
\end{equation}
Conjugating with left translations we get the one-parameter family of
dilations $\varphi_{p,s}:=L_{p}\circ\varphi_{s}\circ L_{p}^{-1}$ with
center at any point $p\in\hh^1$. A set $E\subset\hh^1$ is a \emph{cone
of center $p$} if $\varphi_{p,s}(E)\subset E$ for all $s\in\rr$.

Any isometry of $(\hh^1,g)$ leaving invariant the horizontal
distribution preserves the area of surfaces in $\hh^1$.  Examples of
such isometries are left translations, which act transitively on
$\hh^1$.  The Euclidean rotation of angle $\theta$ about the $t$-axis
given by
\[
(x,y,t)\mapsto r_{\theta}(x,y,t)=(\cos\theta\,x-\sin\theta\,y,
\sin\theta\,x+\cos\theta\,y,t),
\]
is also an area-preserving isometry in $(\hh^1,g)$ since it transforms the
orthonormal basis $\{X,Y,T\}$ at the point $p$ into the orthonormal
basis $\{\cos\theta\,X+\sin\theta\,Y, -\sin\theta\,X+\cos\theta\,Y,
T\}$ at the point $r_{\theta}(p)$.

\section{Examples with one singular line}
\label{sec:line}

Consider the $x$-axis in $\hh^1=\rr^3$parametrized by
$\Ga(v):=(v,0,0)$.  Take a non-increasing continuous function
$\alpha:\rr\to (0,\pi)$.  For every $v\in\rr$, consider two horizontal
halflines $L_{v}^+$, $L_{v}^-$ extending from $\Ga(v)$ with angles
$\alpha(v)$ and $-\alpha(v)$, respectively.  The tangent vectors to
these curves at $\Ga(v)$ are given by
$\cos\alpha(v)\,X_{\Ga(v)}+\sin\alpha(v)\,Y_{\Ga(v)}$ and
$\cos\alpha(v)\,X_{\Ga(v)}-\sin\alpha(v)\,Y_{\Ga(v)}$, respectively.

The parametric equations of this surface are given by
\begin{equation}
\label{eq:parametric}
(v,w)\mapsto
\begin{cases}
(v+w\cos\alpha(v),w\sin\alpha(v),-vw\sin\alpha(v)), & w\ge 0,
\\
(v+|w|\cos\alpha(v),-|w|\sin\alpha(v),v|w|\sin\alpha(v)), & w\le 0,
\end{cases}
\end{equation}
One can eliminate the parameters $v$, $w$ to get the implicit equation
\[
t+xy-y|y|\,\cot\alpha\bigg(\!\!-\frac{t}{y}\bigg)=0.
\]
Letting $\beta:=\cot(\alpha)$, we get that $\beta$ is a continuous
non-decreasing function, and that the surface $\Sg_{\beta}$ defined by
the parametric equations \eqref{eq:parametric} is given by the
implicit equation
\begin{equation}
\label{eq:fbeta}
0=f_{\beta}(x,y,t):=t+xy-y|y|\,\beta\bigg(\!\!-\frac{t}{y}\bigg).
\end{equation}
Observe that, because of the monotonicity condition on $\alpha$, the
projection of relative interiors of the open horizontal halflines to
the $xy$-plane together with the planar $x$-axis $L_{x}$ produce a
partition of the plane.  Since $\Sg_{\beta}$ is the union of the
horizonal lifting of these planar halflines and the $x$-axis to
$\hh^1$, it is the graph of a continuous function
$u_{\beta}:\rr^2\to\rr$.  For $(x,y)\in\rr^2$, the only point in the
intersection of $\Sg_{\beta}$ with the vertical line passing through
$(x,y)$ is precisely $(x,y,u_{\beta}(x,y))$.  Obviously
\begin{equation}
\label{eq:ubeta}
f_{\beta}(x,y,u_{\beta}(x,y))=0.
\end{equation}
For any $(x,y)\in\rr^2$, denote by $\xi_{\beta}(x,y)$ the only value
$v\in\rr$ so that either $\Ga(v)=(x,y,0)$, or $(x,y,u_{\beta}(x,y))$
is contained in one of the two above described halflines leaving
$\Ga(v)$.  Trivially $\xi_{\beta}(x,0)=x$.  Using
\eqref{eq:parametric} one checks that
\begin{equation}
\label{eq:xi}
\xi_{\beta}(x,y)=-\frac{u_{\beta}(x,y)}{y}, \qquad y\neq 0.
\end{equation}
Recalling that $\alpha=\cot^{-1}(\beta)$, we see that the mapping
\[
(v,w)\mapsto 
\begin{cases}
(v+w\cos\alpha(v), w\sin\alpha(v)), &w\ge 0,
\\
(v+|w|\cos\alpha(v), -|w|\sin\alpha(v)), &w\le 0,
\end{cases}
\]
is an homeomorphism of $\rr^2$ whose inverse is given by
\[
(x,y)\mapsto (\xi_{\beta}(x,y), \sgn(y)\,|(x-\xi_{\beta}(x,y), y)|),
\]
where $\sgn(y):=y/|y|$ for $y\neq 0$.  Hence $\xi_{\beta}:\rr^2\to\rr$
is a continuous function.  By \eqref{eq:xi}, the function
$u_{\beta}(x,y)/y$ admits a continuous extension to $\rr^2$.

Let us analyze first the properties of $u_{\beta}$ for regular $\beta$

\begin{lemma}
\label{lem:smooth}
Let $\beta\in C^k(\rr)$, $k\ge 2$, be a non-decreasing function. Then
\begin{enumerate}
\item\label{ck1} $u_{\beta}$ is a $C^k$ function in $\rr^2-L_{x}$,
\item\label{ck2} $u_{\beta}$ is merely $C^{1,1}$ near the $x$-axis
when $\beta\neq 0$,
\item\label{ck3} $u_{\beta}$ is $C^\infty$ in $\xi^{-1}(I)$ when $\beta\equiv
0$ on any open set $I\subset\rr$, and
\item\label{ck4} $\Sg_{\beta}$ is area-minimizing.
\item\label{ck5} The projection of the singular set of $\Sg_{\beta}$
to the $xy$-plane is $L_{x}$.
\end{enumerate}
\end{lemma}

\begin{proof}
Along the proof we shall often drop the subscript $\beta$ for
$f_{\beta}$, $u_{\beta}$, $\xi_{\beta}$ and $\Sg_{\beta}$.

The proof of \ref{ck1} is just an application of the Implicit Function
Theorem since $f_{\beta}$ is a $C^k$ function for $y\neq 0$ when
$\beta$ is $C^k$.

To prove \ref{ck2} we compute the partial derivatives of
$u_{\beta}$ for $y\neq 0$. They are given by
\begin{align}
\label{eq:ux}
(u_{\beta})_{x}(x,y)&=\frac{-y}{1+|y|\,
\beta'\big(\xi_{\beta}(x,y)\big)},
\\
\label{eq:uy}
(u_{\beta})_{y}(x,y)&=\frac{-x+|y|\,\big(2\beta\big(\xi_{\beta}(x,y)
\big)-\beta'\big(\xi_{\beta}(x,y)\big)\,\xi_{\beta}(x,y)\big)}
{1+|y|\,\beta'\big(\xi_{\beta}(x,y)\big)}.
\end{align}
Since $u_{\beta}(x,0)=0$ for all $x\in\rr$ we get
$(u_{\beta})_{x}(x,0)=0$. On the other hand
\[
(u_{\beta})_{y}(x,0)=\lim_{y\to 0}\frac{u_{\beta}(x,y)}{y}=-\lim_{y\to 
0}\xi_{\beta}(x,y)=-\xi_{\beta}(x,0)=-x.
\]
The limits, when $y\to 0$, of \eqref{eq:ux} and \eqref{eq:uy} can be
computed using \eqref{eq:xi}.  We conclude that the first derivatives
of $u_{\beta}$ are continuous functions and so $u_{\beta}$ is a $C^1$
function on $\rr^2$.  To see that $u_{\beta}$ is merely lipschitz, we 
get from \eqref{eq:uy} and \eqref{eq:xi}
\begin{align*}
(u_{\beta})_{yy}(x,0)&=\lim_{y\to 0^\pm}\frac{(u_{\beta})_{y}(x,y)+x}{y}
\\
&=\lim_{y\to 0^\pm}\frac{|y|\,\big(2\beta\big(\xi_{\beta}(x,y)\big)
-\beta'\big(\xi_{\beta}(x,y)\big)\,\xi_{\beta}(x,y)
+x\beta'\big(\xi_{\beta}(x,y)\big)\big)}
{y\,\big(1+|y|\,\beta'\big(\xi_{\beta}(x,y)\big)\big)}
\\
&=\pm 2\beta(x).
\end{align*}
Hence side derivatives exist, but they do not coincide unless
$\beta(x)=0$.

As $u_{\beta}\big|_{\xi^{-1}(I)}=-xy$, \ref{ck3} follows easily .

To prove \ref{ck4} we use a calibration argument. We shall drop the
subscript $\beta$ to simplify the notation. Let $F\subset\hh^1$
such that $F=E$ outside a Euclidean ball $B$ centered at the origin.
Let $H^1:=\{(x,y,t):y\ge 0\}$, $H^2:=\{(x,y,t): y\le 0\}$,
$\Pi:=\{(x,y,t):y=0\}$.  Vertical translations of the horizontal unit
normal $\nu_{E}$, defined outside $\Pi$, provide two
vector fields $U^1$ on $H^1$, and $U^2$ on $H^2$. They are $C^2$ in
the interior of the halfspaces and extend continuously to the boundary
plane $\Pi$. As in the proof of Theorem~5.3 in \cite{rr2}, we see that
\[
\divv U^i=0, \qquad i=1, 2,
\]
in the interior of the halfspaces.  Here $\divv U$ is the riemannian
divergence of the vector field $U$.  Observe that the vector field $Y$
is the riemannian unit normal, and also the horizontal unit normal, to
the plane $\Pi$. We may apply the divergence theorem to get
\begin{align*}
0=\int_{E\cap\intt(H^i)\cap B}\divv U^i
&=\int_{E} \escpr{U^i,\nu_{\intt(H^i)\cap B}}\,d|\ptl(\intt(H^i)\cap B)|
\\
&+\int_{\intt(H^i)\cap B}\escpr{U^i,\nu_{E}}\,d|\ptl E|.
\end{align*}
Let $D:=\Pi\cap \overline{B}$. Then, for
every $p\in D$, we have $\nu_{\intt(H^1)\cap B}=-Y$,
$\nu_{\intt(H^2)\cap B}=Y$, and $U^1=J(v)$, $U^2=J(w)$, where $v-w$ 
is proportional to $Y$, by the construction of $\Sg_{\beta}$. Hence
\[
\escpr{U^1,\nu_{\intt(H^1)\cap B}}+\escpr{U^2,\nu_{\intt(H^2)\cap B}}
=\escpr{v-w,J(Y)}=0, \qquad p\in D.
\]
Adding the above integrals we obtain
\[
0=\sum_{i=1,2}\int_{E}\escpr{U^i,\nu_{B}}\,d|\ptl
B|+\sum_{i=1,2}\int_{B\cap\intt(H^i)}\escpr{U^i,\nu_{E}}\,d|\ptl E|.
\]
We apply the same arguments to the set $F$ and, since $E=F$ on $\ptl
B$ we conclude
\begin{equation}
\label{eq:maincal}
\sum_{i=1,2}\int_{B\cap\intt(H^i)}\escpr{U^i,\nu_{E}}\,d|\ptl
E|=\sum_{i=1,2}\int_{B\cap\intt(H^i)}\escpr{U^i,\nu_{F}}\,d|\ptl F|.
\end{equation}
As $E$ is a subgraph, $|\ptl E|(\Pi)=0$ and so
\[
|\ptl E|(B)=\sum_{i=1,2}\int_{B\cap\intt(H^i)}\escpr{U^i,\nu_{B}}\,d|\ptl E|.
\]
Cauchy-Schwarz inequality and the fact that $|\ptl F|$ is a positive
measure imply
\[
\sum_{i=1,2}\int_{B\cap\intt(H^i)}\escpr{U^i,\nu_{F}}\,d|\ptl F|\le |\ptl F|(B),
\]
which implies \ref{ck4}.

To prove \ref{ck5} simply take into account that the projection of
the singular set of $\Sg_{\beta}$ to the $xy$-plane is composed of
those points $(x,y)$ such that $(u_{\beta})-x-y=(u_{\beta})_{y}+x=0$. 
From \eqref{eq:ux} we get that $(u_{\beta})_{x}-y=0$ if and only if
\[
y\,\big(2+|y|\,\beta'(\xi_{\beta}(x,y))\big)=0,
\]
i.e, when $y=0$.  In this case, from \eqref{eq:uy}, we see that
equation $(u_{\beta})_{y}+x=0$ is trivially satisfied.
\end{proof}


We now prove the general properties of $\Sg_{\beta}$ from Lemma~\ref{lem:smooth}

\begin{proposition}
Let $\beta:\rr\to\rr$ be a continuous non-decreasing function. Let
$u_{\beta}$ be the only solution of equation \eqref{eq:ubeta},
$\Sg_{\beta}$ the graph of $u_{\beta}$, and $E_{\beta}$ the subgraph
of $u_{\beta}$. Then
\begin{enumerate}
\item\label{prop1} $u_{\beta}$ is locally lipschitz in Euclidean sense,
\item\label{prop2} $E_{\beta}$ is a set of locally finite perimeter in
$\hh^1$, and
\item\label{prop3} $\Sg_{\beta}$ is area-minimizing in $\hh^1$.
\end{enumerate}
\end{proposition}

\begin{proof}
Let
\[
\beta_{\eps}(x):=\int_{\rr} \beta(y)\,\eta_{\eps}(x-y)\,dy
\]
the usual convolution, where $\eta$ is a Dirac function and
$\eta_{\eps}(x):=\eta(x/\eps)$, see \cite{eg}.  Then $\beta_{\eps}$ is a
$C^\infty$ non-decreasing function, and $\beta_{\eps}$ converges
uniformly, on compact subsets of $\rr$, to $\beta$.  Let
$u=u_{\beta}$, $u_{\eps}=u_{\beta_{\eps}}$, $f=f_{\beta}$,
$f_{\eps}=f_{\beta_{\eps}}$.

Let $D\subset\rr^2$ be a bounded subset.  To check that $u$ is
lipschitz on $D$ it is enough to prove that the first derivatives of
$u_{\eps}$ are uniformly bounded on $D$.

From \eqref{eq:ubeta} we get
\[
\xi(x,y)+|y|\,\beta\big(\xi(x,y)\big)=x, \qquad y\neq 0.
\]
For $y$ fixed, define the continuous strictly increasing function
\[
\rho_{y}(x):=x+|y|\,\beta(x).
\]
Hence we get
\begin{equation}
\label{eq:rho}
\xi(x,y)=\rho_{y}^{-1}(x).
\end{equation}
We can also define $(\rho_{\eps})_{y}(x):=x+|y|\beta_{\eps}(x)$.  Equation
\eqref{eq:rho} holds replacing $u$, $\beta$ by $u_{\eps}$, $\beta_{\eps}$.

Since $\rho_{y}^{-1}(x)=\xi(x,y)$, we conclude that $\rho_{y}^{-1}$ is
a continuous function that depends continuously on $y$.


Let us estimate 
\[
|(\rho_{\eps})_{y}^{-1}(x)-\rho_{y}^{-1}(x)|.
\]
Let $z_{\eps}:=(\rho_{\eps})_{y}^{-1}(x)$, $z=\rho_{y}^{-1}(x)$. Then 
$x=(\rho_{\eps})_{y}(z_{\eps})=\rho_{y}(z)$ and we have, assuming
$z_{\eps}\ge z$.
\begin{align*}
0&=(\rho_{\eps})_{y}(z_{\eps})-\rho_{y}(z)=
z_{\eps}+|y|\beta_{\eps}(z_{\eps})-\big(z+|y|\beta(z)\big)
\\
&=(z_{\eps}-z)+|y|\,\big(\beta_{\eps}(z_{\eps})-\beta_{\eps}(z)\big)
+|y|\,\big(\beta_{\eps}(z)-\beta(z)\big)
\\
&\ge (z_{\eps}-z)+|y|\,\big(\beta_{\eps}(z)-\beta(z)\big).
\end{align*}
A similar computation can be performed for $z_{\eps}\ge z$. The
consequence is that
\[
|z_{\eps}-z|\le |y|\,|\beta_{\eps}(z)-\beta(z)|,
\]
or, equivalently,
\[
|(\rho_{\eps})_{y}^{-1}(x)-\rho_{y}^{-1}(x)|\le |y|\,
|\beta_{\eps}(\rho_{y}^{-1}(x))-\beta(\rho_{y}^{-1}(x))|.
\]
As $\beta_{\eps}\to\beta$ uniformly on compact subsets of $\rr$, we
have uniform convergence of $(\rho_{\eps})_{y}^{-1}(x)$ to
$\rho_{y}^{-1}(x)$ on compact subsets of $\rr^2$.  This also implies
the uniform convergence of $\xi_{\eps}(x,y)$ to $\xi(x,y)$ on compact
subsets.  Hence also $u_{\eps}(x,y)$ converges uniformly to $u(x,y)$
on compact subsets of $\rr^2$.

From \eqref{eq:ux} and \eqref{eq:uy} we have
\begin{align*}
|(u_{\eps})_{x}(x,y)|&\le |y|, 
\\
|(u_{\eps})_{y}(x,y)|&\le
|x|+2\,|y|\,\big|\beta_{\eps}\big(\xi_{\eps}(x,y)\big)\big|
+\big|\xi_{\eps}(x,y)\big|.
\end{align*}
As $\beta_{\eps}\to\beta$ and $\xi_{\eps}(x,y)\to \xi(x,y)$ uniformly
on compact subsets, we have that the first derivatives of $u_{\eps}$
are uniformly bounded on compact subsets.  Hence $u$ is locally
lipschitz.

%

The subgraph of $u_{\beta}$ is a set of locally finite perimeter in
$\hh^1$ since its boundary is locally lipschitz by \ref{prop1}.  This
follows from \cite{fssc} and proves \ref{prop2}.

To prove \ref{prop3} we use approximation and the calibration argument.
Let $F\subset\hh^1$ so that $F=E$ outside a Euclidean ball $B$
centered at the origin. For the functions $\beta_{\eps}$, consider 
the vector fields $U_{\eps}^i$ obtained by translating vertically the
horizontal unit normal to the surface $\Sg_{\eps}$. We repeat the
arguments on the proof of \ref{ck4} in Lemma~\ref{lem:smooth} to
conclude as in \eqref{eq:maincal} that
\[
\sum_{i=1,2}\int_{B\cap\intt(H^i)}\escpr{U_{\eps}^i,\nu_{E}}\,d|\ptl E|=
\sum_{i=1,2}\int_{B\cap\intt(H^i)}\escpr{U_{\eps}^i,\nu_{F}}\,d|\ptl F|.
\]
Trivially we have
\[
\sum_{i=1,2}\int_{B\cap\intt(H^i)}\escpr{U_{\eps}^i,\nu_{F}}\,d|\ptl F|\le |\ptl
F|(B).
\]
On the other hand, $U_{\eps}^i$ converges uniformly, on compact
subsets, to $U^i$ by Lemma~\ref{lem:nuh}. Passing to the limit when
$\eps\to 0$ and taking into account that $U^i=\nu_{E}$ we conclude
\[
|\ptl E|(B)\le |\ptl F|(B),
\]
as desired.
\end{proof}

\begin{lemma}
\label{lem:nuh}
Let $\beta$ be a continuous non-decreasing function.  Then the
horizontal unit normal of $\Sg_{\beta}$ is given, in
$\{X,Y\}$-coordinates, by
\begin{equation}
\label{eq:nubeta}
\nu_{\beta}(x,y)=\bigg(\frac{1}{(1+\beta^2)^{1/2}},
\frac{-\sgn(y)\,\beta}{(1+\beta^2)^{1/2}}\bigg)\,
\big(\xi_{\beta}(x,y)\big),\qquad y\neq 0.
\end{equation}
Moreover, $\nu_{\beta}$ admits
continuous extensions to $y=0$ from both sides of this line.
\end{lemma}

\begin{proof}
Since $u_{\beta}$ is lipschitz, it is differentiable almost everywhere
on $\rr^2$. On these points,
\[
\nu_{\beta}(x,y)=((u_{\beta})_{x}-y, (u_{\beta})_{y}+x).
\]
The function $-u_{\beta}(x,y)/y$ is constant along the lines
$(x_{0},0)+\lambda\,(1+\beta^2)^{-1/2}(\beta,\pm 1)(x_{0})$, for
$\lambda\ge 0$. Let $y\ge 0$. From \eqref{eq:fbeta} we have
\[
0=-x_{0}+x-y\,\beta(x_{0}).
\]
Let $v:=(1+\beta^2)^{-1/2}(\beta,1)(x_{0})$. Then
$v(-u_{\beta}(x,y)/y)=0$. Hence for almost every point on almost every
line, we have 
\[
\beta(x_{0})\,(u_{\beta})_{x}+(u_{\beta})_{y}=-x_{0}.
\]
Hence we have
\[
(u_{\beta})_{y}+x=-x_{0}-\beta(x_{0})\,(u_{\beta})_{x}+x_{0}+y\,\beta(x_{0})
=\beta(x_{0})\,(-(u_{\beta})_{x}+y).
\]
We conclude that the horizontal unit normal is proportional to
$(1,-\beta)$, which implies \eqref{eq:nubeta}. The case $y\le 0$ is
handled similarly.
\end{proof}


\begin{example}
Taking $\beta(x):=x$ we get
\[
u_{\beta}(x,y)=-\frac{xy}{1+|y|},
\]
which is a Euclidean $C^{1,1}$ graph.

Another family of interesting examples are the minimal cones obtained 
by taking the constant function $\beta(x):=\beta_{0}$. In this case we get
\[
u_{\beta}(x,y)=-xy+\beta_{0}\,y|y|.
\]
In this case $\Sg_{\beta}$ is a $C^{1,1}$ surface which is invariant
by the dilations centered at any point of the singular line.

Take now
\[
\beta(x):=\begin{cases}
0, & x\le 0, \\
x, & x\ge 0.
\end{cases}
\]
In this case we obtain the graph
\[
u_{\beta}(x,y):=\begin{cases}
-xy, & x\le 0, \\
-\frac{xy}{1+|y|}, & x\ge 0,
\end{cases}
\]
which is simply locally Lipschitz.

This example was mentioned to me by Scott Pauls.  Consider now a
continuous nondecreasing function $\beta:\rr\to\rr$, constant outside
the Cantor set $C\subset [0,1]$ with $\beta(0)=0$, $\beta(1)=1$.  Then
the associated surface $\Sg_{\beta}$ is an area-minimizing surface in
$\hh^1$.
\end{example}

%
%
%
%
%

\section{Examples with several singular halflines meeting at a point}
\label{sec:halflines}

Let $\alpha_{1}^0, \ldots,\al_{k}^0$, be a family of positive angles
so that
\[
\sum_{i=1}^k\alpha_{i}^0=\pi.
\]
Let $r_{\beta}$ be the rotation of angle $\beta$ around the origin in
$\rr^2$.  Consider a family of closed halflines $L_{i}\subset\rr^2$,
$i\in\zz_{k}$, extending from the origin, so that
$r_{\alpha_{i}^0+\alpha_{i+1}^0}(L_{i})=L_{i+1}$.  Finally, define
$R_{i}:=r_{\alpha_{i}^0}(L_{i})$.  (An alternative way of defining this
configuration is to start from a family of counter-clockwise oriented
halflines $R_{i}\subset\rr^2$, $i\in\zz_{k}$, choosing $L_{i}$,
$i\in\zz_{k}$, as the bisector of the angle determined by $R_{i-1}$
and $R_{i}$, and defining $\alpha_{i}^0$ as the angle between $L_{i}$
and $R_{i}$). Define $W_{i}$ as the closed wedge, containing $L_{i}$,
bordered by $R_{i-1}$ and $R_{i}$.

\begin{figure}[h]
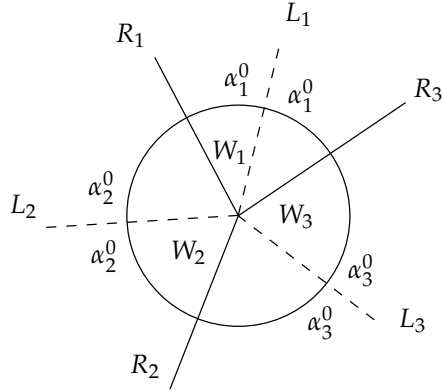
\caption{The initial configuration with three halflines $L_{1}$,
$L_{2}$, $L_{3}$.}
\end{figure}

For every $i\in\zz_{k}$, let $\alpha_{i}:[0,\infty)\to (0,\pi)$ be a
continuous nonincreasing function so that $\al_{i}(0)=\al_{i}^0$, and
define, as in the previous section, $\beta_{i}:=\cot(\al_{i})$.  Let
$v_{i}\in\sph^1$, $i\in\zz_{k}$, be such that $L_{i}=\{sv_{i}: s\ge
0\}$.  For every $i\in\zz_{k}$ and $s\ge 0$, we take the two closed
halflines $L_{s,i}^\pm$ in $\rr^2$ extending from the point $sv_{i}$
with tangent vectors $(\cos\al_{i}(s),\pm\sin\al_{i}(s))$.  In this
way we cover all of $\rr^2$. We shall define
$\alpha:=(\alpha_{1},\ldots,\al_{k})$.

Lift $L_{1},\ldots,L_{k}$ to horizontal halflines
$L_{1}',\ldots,L_{k}'$ in $\hh^1$ from the origin, and $L_{s,i}^\pm$
to horizontal halflines in $\hh^1$ extending from the unique point in
$L_{i}'$ projecting onto $sv_{i}$.  In this way we obtain a continuous
function $u_{\alpha}:\rr^2\to\rr$.  The graph $\Sg_{\al}$ of $u_{\al}$
is a topological surface in $\hh^1$.

Obviously the angle functions $\alpha_{i}(s)$ can be extended
continuously and preserving the monotonicity, to an angle function
$\wt{\al}_{i}:\wt{L}_{i}\to (0,\pi)$, where $\wt{L}_{i}$ is the
straight line containing the halfline $L_{i}$. The graph of $u_{\al}$ 
restricted to $W_{i}$ coincides with the Euclidean locally lipschitz
area-minimizing surface $u_{\wt{\beta_{i}}}$, for
$\wt{\beta}_{i}:=\cot{\wt{\alpha}_{i}}$, constructed in the previous
section. So the examples in this section can be seen as pieces of the 
examples of the previous one patched together. 

\begin{theorem}
Under the above conditions
\begin{enumerate}
\item The function $u_{\al}$ is locally lipschitz in the Euclidean 
sense.
\item The surface $\Sg_{\al}$ is area-minimizing.
\end{enumerate}
\end{theorem}

\begin{proof}
It is immediate that $u_{\alpha}$ is a graph which is locally
lipschitz in Euclidean sense: choose a disk $D\subset\rr^2$.  Let $p$,
$q\in D$.  Assume first that $(p,q)$ intersects the halflines
$R_{1},\ldots,R_{k}$ transversally at the points $x_{1},\ldots,
x_{n}$.  Then $[p,x_{1}]$, $[x_{1},x_{2}], \ldots , [x_{n},p]$ are
contained in wedges and hence
\begin{align*}
|u_{\al}(p)-u_{\al}(q)|&\le |u_{\al}(p)-u_{\al}(x_{1})|+ \cdots
+|u_{\al}(x_{n})-u_{\al}(q)|
\\
&C\,\big(|p-x_{1}|+\cdots |x_{n}-q|\big)=C\,|p-q|,
\end{align*}
where $C$ is the supremum of the Lipschitz constants of
$u_{\wt{\beta}_{i}}$ restricted to $D$.  The general case is then
obtained by approximating $p$ and $q$ by points in the condition of
the assumption.

To prove that $u_{\al}$ is area minimizing we first approximate
$\alpha_{i}$ by smooth angle functions $(\alpha_{i})_{\eps}$ with
$(\alpha_{i})_{\eps}(0)=\alpha_{i}(0)$. In this way we obtain a
calibrating vector field which is continuous along the vertical planes
passing through $R_{i}$ by Lemma~\ref{lem:nuh}. This allows us to
apply the calibration argument to prove the area-minimizing property
of $\Sg_{\alpha}$.
\end{proof}

\begin{example}[Minimizing cones]
Let $\al_{i}(s)=\al_{i}^{0}$ be a constant for all $i$.  Then the
subgraph of $\Sg_{\alpha}$ is a minimizing cone with center at $0$.
Restricted to the interior of the wedges $W_{i}$, the surface
$\Sg_{\alpha}$ is $C^{1,1}$.  An easy computation shows that, taking
$\beta(s):=\beta_{0}$ in the construction of the first section, the
Riemannian normal to $\Sg_{\beta}$ along the halflines
$\beta_{0}|y|=x$, $x\ge 0$ (that make angle $\pm\cot^{-1}(\beta_{0})$
with the positive $x$-axis) is given by
\[
N=\frac{-2y\,X+2\beta_{0}|y|\,Y-T}{\sqrt{1+4y^2+4\beta_{0}^2y^2}}=
\frac{-2y\,X+2x\,Y-T}{\sqrt{1+4x^2+4y^2}}.
\]
This vector field is invariant by rotations around the vertical axis. 
Hence in our construction, the normal vector field to $\Sg_{\alpha}$
is continuous. It is straightforward to show that it is locally
lipschitz in Euclidean sense.
\end{example}

\begin{example}[Area-minimizing surfaces with a singular halfline]
These examples are inspired by \cite[Example~7.2]{chy1}.  We consider
a halfline $L$ extending from the origin, and an angle function
$\alpha:L\to (0,\pi)$ continuous and nonincreasing as a function of
the distance to the origin.  We consider the union of the halflines
$L_{\alpha(q)}^+$, $L_{\alpha(q)}^-$ extending from $q\in L$ with
angles $\alpha(q)$, $-\alpha(q)$, respectively.  We patch the
area-minimizing surface defined by $\alpha$ in the wedge delimited by
the halflines $L^+_{0}$, $L^-_{0}$, with the plane $t=0$.  In this way
we get an entire area-minimizing $t$-graph, with lipschitz regularity.
In case the angle function $\alpha$ is constant, we get an
area-minimizing cone with center $0$, which is defined by the equation
\[
u(x,y):=\begin{cases}
-xy+\beta_{0}\,y|y|, &-xy+\beta_{0}\,y|y|\ge 0,
\\
0, &-xy+\beta_{0}\,y|y|\le 0.
\end{cases}
\]
This surface is composed of two smooth pieces patched together along 
the halflines $x=\beta_{0}|y|$.
\end{example}

\providecommand{\bysame}{\leavevmode\hbox to3em{\hrulefill}\thinspace}
\providecommand{\MR}{\relax\ifhmode\unskip\space\fi MR }
\providecommand{\MRhref}[2]{%
  \href{http://www.ams.org/mathscinet-getitem?mr=#1}{#2}
}
\providecommand{\href}[2]{#2}


\end{document}